\documentclass[amscd,amssymb,verbatim,12pt]{amsart}
\usepackage{amsmath, amsthm, amssymb}
\usepackage{epsfig}
\usepackage[iso-8859-7]{inputenc}
\usepackage{enumerate}
\usepackage{mathrsfs}
\input{epsf.tex}
\newtheorem{lemma}{Lemma}
\newtheorem{sub}{Lemma}[lemma]
\newtheorem{theorem}{Theorem}
\newtheorem{defin}{Definition}
\newtheorem*{Alex}{Alexander's Lemma (for the plane)}
\newtheorem*{Alex-sc}{Alexander's Lemma (for simply connected spaces)}

\newtheorem{question}{Question}
\newcommand{\diam}{\textnormal{diam}}
\newcommand{\al}{\alpha}
\newcommand{\sm}{\smallsetminus}
\newcommand{\ep}{\varepsilon}

\newcommand{\os}{\caption{}}

\begin{document}

\title[\fontsize{7}{0}\selectfont Simply connected homogeneous continua are not separated by arcs]{Simply connected homogeneous continua are not separated by arcs}

\author[\fontsize{7}{0}\selectfont Myrto Kallipoliti and Panos Papasoglu]
{Myrto Kallipoliti and Panos Papasoglu}


\email [Myrto Kallipoliti]{mirtok@math.uoa.gr} \email [Panos
Papasoglu]{panos@math.uoa.gr}

\address
[Myrto Kallipoliti] {Mathematics Department, University of Athens,
Athens 157 84, Greece }
\address
[Panos Papasoglu] {Mathematics Department, University of Athens,
Athens 157 84, Greece }

\thanks{Research co-funded by European Social Fund and National Resources EPEAEK II-Pythagoras}

\date{\today}


\begin{abstract}
We show that locally connected, simply connected homogeneous
continua are not separated by arcs.
We ask several questions about homogeneous
continua which are inspired by analogous questions in geometric
group theory.
\end{abstract}

\maketitle

\section{Introduction}

In this paper we prove a theorem about homogeneous continua
inspired by a result about finitely presented groups (\cite {P}).

\begin{theorem}
\label{continua} Let $X$ be a locally connected, simply connected,
homogeneous continuum. Then no arc separates $X$.
\end{theorem}

We recall that an arc in $X$ is the image of a 1-1 continuous map
$\alpha :[0,1]\to X$. We say that an arc $\alpha$ separates $X$ if
$X\smallsetminus\alpha$ has at least two connected components. We
say that $X$ is simply connected if it is path connected and every
continuous map $f:S^1=\partial D^2\to X$ can be extended to a
continuous $\bar f:D^2\to X$ (where $D^2$ is the 2-disc and $S^1$
its boundary circle). The proof of Theorem \ref{continua} relies
on Alexander's lemma for the plane (see \cite {N}) and our
generalization of this lemma to simply connected spaces (see sec.
2).

There is a family air between continua theory and group theory.
This became apparent after Gromov's theory of hyperbolic groups (\cite {Gr1}).
Gromov defines a boundary for a hyperbolic group which is a
continuum on which the group acts by a `convergence action'. The
classic `cyclic elements' decomposition theory of Whyburn was
extended recently by Bowditch (\cite {Bo}) in this context and it
gave deep results in group theory. `Asymptotic topology',
introduced by Gromov (\cite {Gr2}) and developed further by
Dranishnikov (\cite {Dr}), shows that the analogy goes beyond the
realm of hyperbolic groups. The `philosophy' of this is that
topological questions that make sense for continua can be
translated to `asymptotic topology' questions which make sense for
groups (see \cite {Dr} for a dictionary between topology and
asymptotic topology).

One wonders whether Theorem \ref{continua} holds in fact for all
locally connected, homogeneous continua of dimension bigger than
1:

\begin{question}
Let $X$ be a locally connected, homogeneous continuum of dimension
2. Is it true that no arc separates $X$?
\end{question}

We note that by a result of Krupski (\cite{Kr}), homogeneous
continua are Cantor manifolds. It follows that no arc separates a
homogeneous continuum of dimension bigger than 2.

 Krupski
and Patkowska (\cite {K-P}) have shown that a similar property
(the disjoint arcs property) holds for all locally connected
homogeneous continua of dimension bigger than 1 which are not
2-manifolds.

We remark that by \cite{Kr} if a Cantor set separates a
homogeneous continuum $X$ then $\dim X=1$. So question 1 is
equivalent to the following question: Is it true that if no Cantor
set separates a locally connected, homogeneous continuum $X$, then
no arc separates $X$? Restated in this way the question makes
sense also for boundaries of hyperbolic groups. In fact a similar
question can be formulated for finitely generated groups too (see
\cite {P}).

Not much is known about locally connected, simply connected
homogeneous continua. One motivation to study them is the analogy
with finitely presented groups. Another reason is that one could
hope for a classification of such continua in dimension 2:

\begin{question}
Are the 2-sphere and the universal Menger compactum of dimension 2
the only locally connected, simply connected, homogeneous continua
of dimension 2?
\end{question}

We recall that $S^1$ and the universal Menger curve are the only
locally connected homogeneous continua of dimension 1 (\cite {A}).
Prajs (\cite {Pr1}, question 2) asks whether $S^2$ is the only
simply connected homogeneous continuum of dimension 2 that embeds
in $\mathbb{R }^3$.

A related question about locally connected, simply connected
continua that makes sense also for finitely presented groups is
the following:

\begin{question}
Let $X$ be a locally connected, simply connected homogeneous
continuum which is not a single point. Does $X$ contain a disc?
\end{question}
We remark that in the group theoretic setting the answer is
affirmative for hyperbolic groups (\cite {B-K}). By a result of
Prajs (\cite {Pr2}) a positive answer to this would imply that
$S^2$ is the only locally connected, simply connected, homogeneous
continuum of dimension 2 that embeds in $\mathbb{R}^3$.

We refer to Prajs' list of problems (\cite {Pr1}) for more
questions on homogeneous continua.

\section{Preliminaries}

\begin{defin}
Let $X$ be a metric space. A path $p$ is a continuous
map $p:[0,1]\to X$. A simple path or an arc $\al$, is a continuous
and $1-1$ map $\al:[0,1]\to X$. We will identify an arc with its
image.

\

\noindent For a path $p$ we denote by $\partial p$ the set of its endpoints,
i.e. $\partial p=\{p(0),p(1)\}$.

\

\noindent An arc $\al$ separates $X$ if $X\smallsetminus \al$ has
at least two connected components. If $x,y\in X$ we say that an
arc $\al$ separates $x$ from $y$ if $\al$ separates $X$ and $x,y$
belong to distinct components of $X\smallsetminus \al$.

\end{defin}

\

\begin{defin}
Let $\al$ be an arc of $X$. On $\al$ we define an order
$<_{\al}$ as follows:
 If $x=\al(x'),\,y=\al(y')$ then $x<_{\al}y$ if and only
if $x'<y'$.

\

\noindent We denote by ${\left[x,y\right]}_{\al}$ the set of all
$t\in {\al}$ such that $x\leq t\leq y$. Similarly we define
${(x,y)}_{\al},\,[x,y)_{\al}$ and $(x,y]_{\al}$. When there is no
ambiguity we write $\left[x,y\right]$ instead of
${\left[x,y\right]}_{\al}$ and $x<y$ instead of $x<_{\al}y$.
Finally, if $t\in [0,1]$ we denote by $x+t$ the point $\al(x'+t)$
(where $x=\al(x')$).
\end{defin}

\

We recall Alexander's lemma from plane topology 
(see Theorem 9.2, p.112 of \cite {N}).

\begin{Alex}
Let $K_1,K_2$ be closed sets on the plane such that either
$K_1\cap K_2=\emptyset $ or $K_1\cap K_2$ is connected and at
least one of $K_1,K_2$ is bounded. Let $x,y\in
\mathbb{R}^2\smallsetminus (K_1\cup K_2)$. If there is a path
joining $x,y$ in $\mathbb{R}^2\smallsetminus K_1$ and a path
joining $x,y$ in $\mathbb{R}^2\smallsetminus K_2$ then there is a
path joining $x,y$ in $\mathbb{R}^2\smallsetminus(K_1\cup K_2)$.
\end{Alex}

It is easy to see that Alexander's lemma also holds for the closed
disc $D^2$ in the case that $K_1\cap K_2=\emptyset $. In fact this
implies that this lemma holds in general for every simply
connected space. In particular we have the following:

\begin{Alex-sc}
Let $X$ be a simply connected space, $K_1,K_2$ disjoint closed
subsets of $X$ and let $x,y\in X\smallsetminus (K_1\cup K_2)$. If
there is a path joining $x,y$ in $X\smallsetminus K_1$ and a path
joining $x,y$ in $X\smallsetminus K_2$ then there is a path
joining $x,y$ in $X\smallsetminus(K_1\cup K_2)$.
\end{Alex-sc}

\begin{proof}
Let $p_1,p_2$ be paths joining $x,y$ such that $p_1\cap
K_1=p_2\cap K_2=\emptyset$. We consider the closed path $p_1\cup
p_2$ and let $f:S^1=\partial D^2\to X$ be a parametrization of
this path. Since $X$ is simply connected, $f$ can be extended to a
map $F:D^2\to X$. Then $F^{-1}(K_1)$ and $F^{-1}(K_2)$ are
disjoint, closed subsets of $D^2$. Clearly, neither $F^{-1}(K_1)$
nor $F^{-1}(K_2)$ separates $x,y$, therefore, using Alexander's
lemma for the closed disc, we have that there is a path $p$ that
joins $x,y$ without meeting $F^{-1}(K_2)\cup F^{-1}(K_2)$. This
implies that $F(p)$ is a path from $x$ to $y$ that does not meet
$K_1\cup K_2$.
\end{proof}

\

For the rest of this paper we assume that $X$ is a simply connected,
locally connected continuum.

\

\begin{lemma}
\label{A.0} Let $O$ be a connected open subset of $X,\ K$ be a
connected component of $\partial O$ and let $x,y\in O$ such that
$d(x,K)<\ep$ and $ d(y,K)<\ep$. Then there is a path $p$ in $O$
connecting $x$ to $y$ such that $p$ is contained in the
$\ep$-neighborhood of $\partial O$.
\end{lemma}

\begin{proof}
Let $U$ be the union of the open balls $B_{\ep}(t)$ with center
$t\in \partial O$ and radius $\ep $. Let $V$ be the connected
component of $U$ containing $K$. Clearly $x,y\in V$ so there is a
path in $X$ joining them that does not intersect $\partial V$. On
the other hand $x,y\in O$ so there is a path in $X$ joining them
that does not intersect $\partial O$.

Since $\partial O \cap \partial V=\emptyset$ and $\partial O,
\partial V$ are closed, applying Alexander's lemma for the simply
connected space $X$, we have that $p$ is a path lying in $X$
joining $x,y$ that  intersects neither $\partial O$ nor $\partial
V$. Clearly $p$ is contained in $O$ and lies in the $\ep$-neighborhood of $\partial O$.
\end{proof}

\begin{lemma}
\label{simplycon} Let $\al$ be an arc that separates $X$ and let
$C$ be a connected component of $X\sm \al$. Then
$\overline C$ is simply connected and $\partial C$ is connected.
\end{lemma}

\begin{proof}
Let $f:S^1 = \partial D^2\to\overline C$. We will show that this
map can be extended to a map $\widehat{f}: D^2 \rightarrow
\overline C$.

$X$ is simply connected, so there is a map $F:D^2 \rightarrow X$
such that $F|_{S^1}=f$. Furthermore, $X\sm \overline C$
is an open set, therefore $\partial F^{-1}(X\sm
\overline C)\cap F^{-1}(X\sm \overline C)=\emptyset$
and since $F$ is a continuous extension of $f$ it follows that
$F(\partial F^{-1}(X\sm \overline C)) \subset \al$
(where by $\partial F^{-1}(X\sm \overline C)$ we denote
the boundary of $F^{-1}(X\sm \overline C)$ in $X$).

Let $f': \partial F^{-1}(X\smallsetminus \overline C) \rightarrow
\al$ be the restriction of $F$ in $\partial F^{-1}(X\smallsetminus
\overline C)$. Then, applying Tietze's extension theorem, we
obtain an extension for $f'$:
 $$F':\overline{F^{-1}(X\smallsetminus \overline C)} \rightarrow \al.$$

\noindent Finally we define $\widehat{f} : D^2 \rightarrow
\overline C$ as follows:

\[\widehat{f}(x)=\left\{\begin{array}{ll}

f(x) & \mbox{if $x\in \partial D^2$}, \\
F(x) & \mbox{if $x\in D^2\smallsetminus F^{-1}(X\smallsetminus \overline C)$}, \\
F'(x) & \mbox{if $x\in F^{-1}(X\smallsetminus \overline C)$}.

\end{array}
\right.\]

\noindent This shows that $\overline C$ is simply connected.

\

Suppose now that $S=\partial C$ is not connected and let $p$ be
a path that joins two different components of $S$, such that
if $a,b$ are the endpoints of $p$, then $(p\sm\{a,b\})\cap S=\emptyset$.
Let $x\in p\sm \{a,b\}$ and $y\in {(a,b)}_{\al}\sm S$.
 (Figure \ref{fig:conn}).

\begin{figure}[h]
\begin{center}
\includegraphics[width=3.6in]{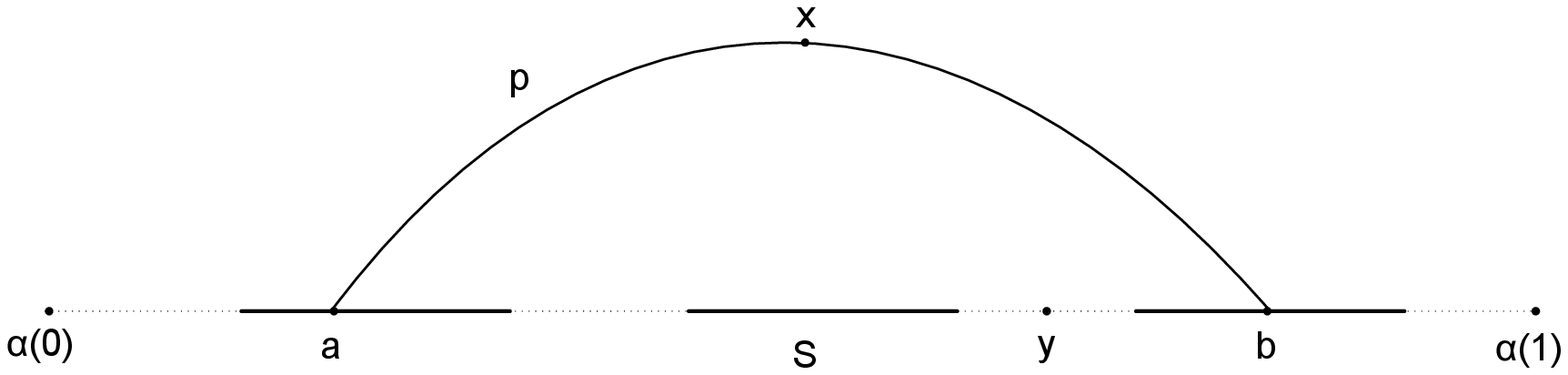}
\end{center}
\os
\label{fig:conn}
\end{figure}

We set $K_1={\left[\al(0),y\right]}_{\al}\cap S$
and $K_2=S\sm K_1$.
It is clear that $K_1,K_2$ are disjoint closed subsets
of $X$ and that neither $K_1$ nor $K_2$ separates $x$ from $y$.
Then Alexander's lemma (for the simply connected space $X$) implies that
there is a path joining $x,y$ in $X\sm (K_1\cup K_2)=X\sm S$, a contradiction.
\end{proof}

\begin{lemma}
\label{A.0'} Let $\al$ be an arc that separates $X,\ x,y\in \al$
and $\ep>0$ with $\ep<d(x,y)$. Then for every connected component
$C$ of $X\smallsetminus \al$ such that $x,y\in\partial C$ there
are points $x',y'\in C$ with $d(x,x'),d(y,y')<\ep$ and a path
$p\in C$ that joins $x',y'$ and is contained in the
$\ep$-neighborhood of ${\left[x,y\right]}_{\al}$.
\end{lemma}

\begin{proof}
Let $B_{\frac{\ep}{2}}(x)$ and $B_{\ep}(x)$ be balls of center $x$
and radius $\frac{\ep}{2}$ and $\ep$ respectively. We consider the
connected components of $\al\sm \overset {\circ
}{B}_{\frac{\ep}{2}}(x)$ and we restrict to those that are not
contained in $B_{\ep}(x)$ (here we denote by $\overset {\circ
}{B}_{\frac{\ep}{2}}(x)$ the open ball). It is clear that there
are finitely many such components, so we denote them by
$I_1,I_2,\dots,I_n$. Let $\delta_1<\min\{d(I_i,I_j)\}$ for every
$i,j=1,2,\dots,n,\ i\neq j$. Similarly, let $J_1,J_2,\dots,J_m$ be
the connected components of $\al\sm\overset {\circ}B_{\frac{\ep}{2}}(y)$
that are not contained in $B_{\ep}(y)$ and
let $\delta_2<\min\{d(J_i,J_j)\}$ for every $i,j=1,2,\dots,m,\
i\neq j$. Let $\delta'<\min\{\delta_1,\delta_2,\frac{\ep}{2}\}$.

From Lemma \ref{simplycon}, we have that $\overline{C}$ is
simply connected, therefore Lemma \ref{A.0}, for
$\delta=\frac{\delta'}{4}$, implies that there is a path $q\in C$
that joins a point of $B_{\ep}(x)$ with a point of $B_{\ep}(y)$ and
lies in the $\delta$-neighborhood of $\al$ (Figure \ref{fig:path}).
We will show that there is a subpath of $q$ that has the required
properties.

\begin{figure}[h]
\begin{center}
\includegraphics[width=3.5in]{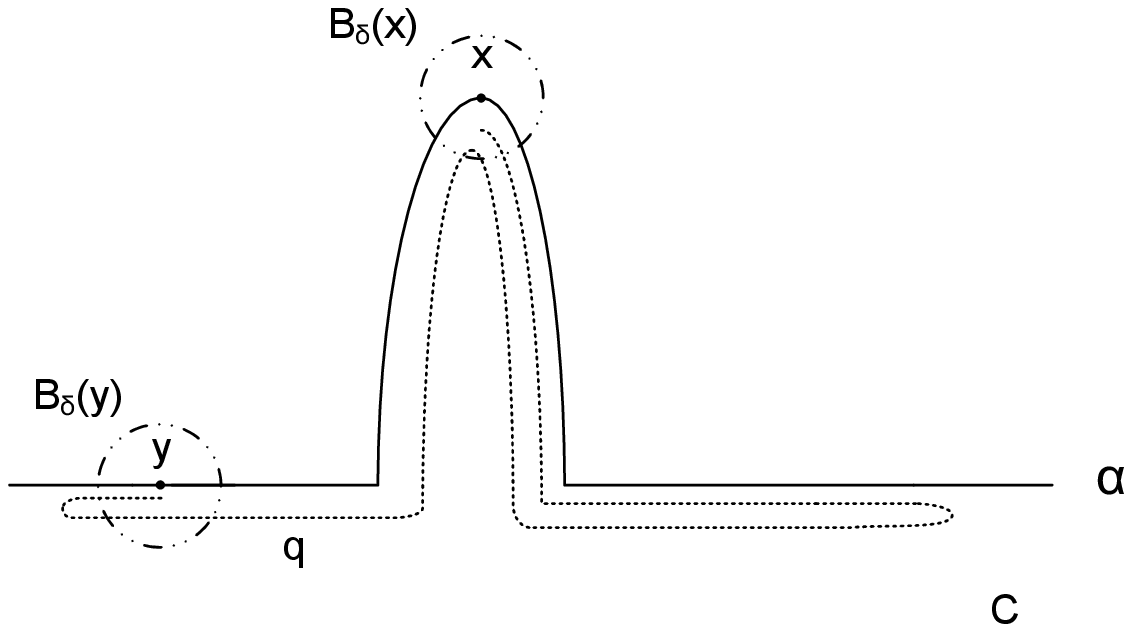}
\end{center}
\os

\label{fig:path}
\end{figure}

We assume that none of the $I_i,J_j,\ i=1,2,\dots,n,\, j=1,2,\dots,m$
contain ${\left[x,y\right]}_{\al}$, since otherwise we are done.
Thus, without loss of generality, let $I_1$ be the connected component
of $\al\sm\overset {\circ} B_{\frac{\ep}{2}}(x)$ that is contained in
${\left[x,y\right]}_{\al}$.
We denote by $N_{\delta}({\left[x,y\right]}_{\al})$
the open $\delta$-neighborhood of ${\left[x,y\right]}_{\al}$.
Suppose that there is a connected component $I=[a,b]_q$ of
$q\smallsetminus  N_{\delta}({\left[x,y\right]}_{\al})$ with $a\in
B_{\delta}(x)$ and $b\notin B_{\delta}(x)$. Then there is an $r>0$
such that ${(b-r,b)}_q\notin N_{\delta}(\al)$. Indeed, if not then
for every $r>0$ there is an $I_i\neq I_1$ such that
${(b-r,b)}_q\in N_{\delta}(I_i)$. Thus $d(I_i,b)\leq \delta$. But
$d(I_1,I_i)\leq d(I_i,b)+d(b,I_1)\leq
\delta+\delta=2\delta=\frac{\delta'}{2}<\delta'$, a contradiction.
So there is an $r>0$ such that for every $i\neq 1$ we have
${(b-r,b)}_q\notin N_{\delta}(I_i)$, therefore ${(b-r,b)}_q\notin
N_{\delta}(\al)$, which is not possible.
This contradiction proves the lemma.
\end{proof}

\

\section{Proof of Theorem \ref{continua}}

\noindent We will prove the theorem by contradiction.

\

\noindent\textbf{Remark: }\normalfont Since $X$ is locally
connected and compact, it follows that every open connected subset
of $X$ is path connected (see Theorem 3.15, p.116 of \cite{H-Y}).
In particular the closure of every component of $X\sm \al$ is
path connected.

\

\begin{defin}
Let $\al_1,\al_2$ be arcs that separate $X$. We say that $\al_1$
crosses $\al_2$ at $x\in (\al_1\sm\partial\al_1)\cap
(\al_2\sm\partial\al_2)$ if for any neighborhood of $x$ in
$\al_2$, $(x-\ep,x+\ep)_{\al_2}$, there are $a,b\in
(x-\ep,x+\ep)_{\al_2}$ separated by $\al_1$. More generally, if
$[x_1,x_2]$ is a connected component of $\al_1\cap \al_2$,
which is contained
in  $(\al_1\sm\partial\al_1)\cap
(\al_2\sm\partial\al_2)$, we say
that $\al_1$ crosses $\al_2$ at $[x_1,x_2]$ if for any
neighborhood of $[x_1,x_2]$ in $\al_2$,
$(x_1-\ep,x_2+\ep)_{\al_2}$, there are $a,b\in
(x_1-\ep,x_2+\ep)_{\al_2}$ separated by $\al_1$. In this case,
the endpoints $x_1,x_2$ are also called cross points of $\al_1,\al_2$.

If $I_1\subset \al_1$, $I_2\subset \al_2$ are intervals of
$\al_1,\al_2$ containing $x$ in their interior, we say that
$I_1,I_2$ cross at $x$. Similarly we define what it means for two
intervals to cross at a common subarc. We call $x$ (respectively
$[x_1,x_2]$) a cross-point (respectively cross-interval) of
$\al_1,\al_2$. We say that $I_1,I_2$ cross if they cross at some
point $x$ or at some interval $[x_1,x_2]$.
\end{defin}

\begin{figure}[h]
\begin{center}
\includegraphics[width=3in]{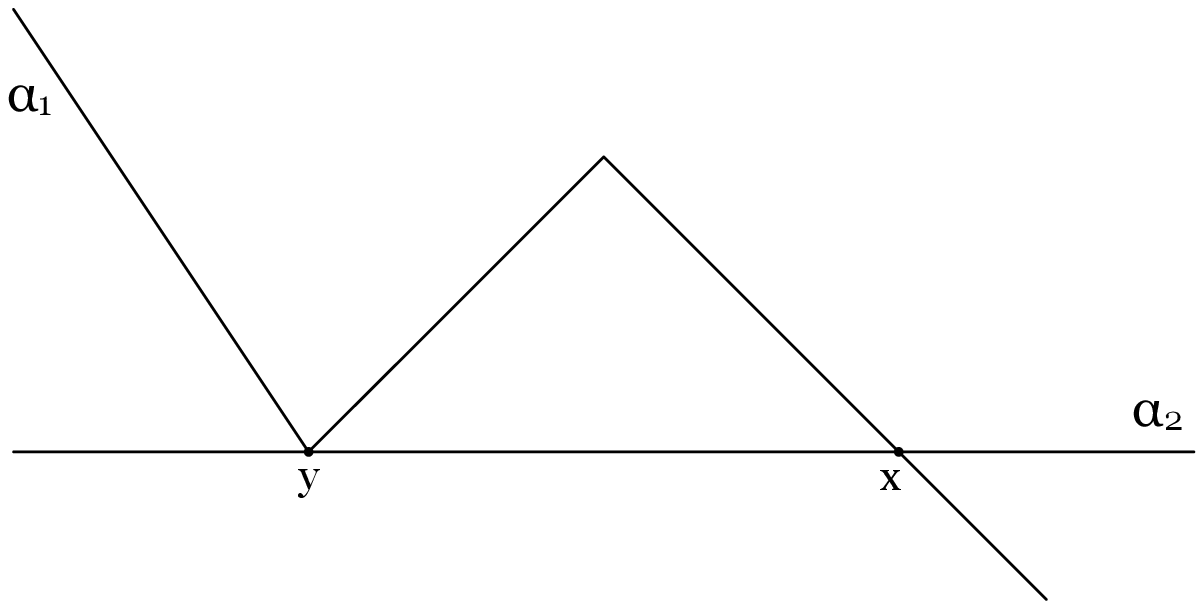}
\end{center}
\os
\label{fig:defcross}
\end{figure}

For example in  Figure \ref{fig:defcross},
$x$ is a cross point of $\al_1,\al_2$,
while $y$ is an intersection point of $\al_1,\al_2$ which is not
a cross point.

\

\begin{lemma}
\label{A.1}
There is an arc that separates $X$ in exactly two
components.
\end{lemma}

\begin{proof}
Suppose that this is not the case, so let $\al$ be an arc that
separates $X$ in more than two components. Since $X$ has no cut
points there are two connected components of $X\sm\al$, say
$C_1,C_2$, such that $\beta=\partial C_1\cap \partial
C_2\neq\emptyset$ is a subarc (which is not a point) of $\al$.
Clearly $\beta$ separates $X$. To simplify notation we denote by
$C_1,C_2$ the components of $X\smallsetminus\beta$ that satisfy
$\partial C_1=\partial C_2=\beta$. Let $C_3$ be another component
of $X\smallsetminus\beta$. By Lemma \ref{simplycon} we have that
$\partial C_3$ is connected, so $\partial C_3=\gamma$ is a subarc
of $\beta$, which separates $X$.

\begin{sub}
\label{subl}
$\gamma$ cannot be crossed by any other separating arc of $X$.
\end{sub}

\begin{proof}
Suppose that there is an arc $\gamma'$ that separates $X$ and
crosses $\gamma$ at $t$. Then there are $x,y\in\gamma,\ x<t<y$
that are separated by $\gamma'$. Let $Y=\overline C_1\cup\overline
C_2$. Since $\overline C_3$ is path connected, it follows that
$\gamma'\not\subset Y$. We denote by $C_x$ the connected component
of $Y\sm \gamma'$ that contains $x$. As in Lemma \ref{simplycon},
we may show that $Y$ is simply connected. We show then that
Alexander's lemma for $Y$ implies that $\partial C_x$ has a
connected component that separates $x,y$ in $Y$.

This can be achieved as follows: Let $C_y$ be the connected
component of $Y\sm \overline{C}_x$ that contains $y$. It is clear
that $\partial C_y\subseteq \partial C_x$. Then no proper closed
subset of $\partial C_y$ separates $x$ from $y$.
Indeed, suppose that there is a closed $K\subset \partial C_y$
that separates $x$
from $y$ and let $z\in \partial C_y\sm K$. Let $U$ be an open
neighborhood of $z$ such that $U\cap K=\emptyset$. It is obvious
that $U$ intersects every component of $Y\sm \partial C_y$,
therefore, there are paths $q_1,q_2\in Y\sm \partial C_y$ that
join $x,y$ with points $x',y'\in U$ respectively. However $U$ is
path connected and since $U\cap K=\emptyset$, it follows that
there is also a path $q\in U$ that joins $x'$ with $y'$. Thus
$p_1\cup q \cup p_2$ is a path joining $x,y$ without meeting $K$,
a contradiction. Therefore, $I=\partial C_y$ is 
connected and separates $x$ from $y$.

We note now that $I$ does not cross $\gamma$. Indeed,
suppose that there is an $a\in I\sm\partial I$ in which $\gamma '$
crosses $\gamma$. Let $V\subset X$ be sufficiently small
neighborhood of $a$  such that $(\gamma'\sm I)\cap V=\emptyset$.
We denote by $J$ the connected component of $\gamma\cap V$ that
contains the point $a$. Then we can pick points $x',y'\in J$ with
$x'<a<y'$ in $\gamma$ that are separated by $\gamma'$. Let
$N_{x'},N_{y'}$ be connected neighborhoods of $x'$ and $y'$
respectively, such that $N_{x'}, N_{y'}\subset V$. Applying now
Lemma \ref{A.0'} for the component $C_3$ and for $\ep
<\min\{\diam(N_{x'}\cap C_3),\diam(N_{y'}\cap C_3)\}$, we have
that every point of $N_{x'}\cap C_3$ can be joined with every
point of $N_{y'}\cap C_3$ by a path in $C_3$ which lies in the
$\ep$-neighborhood of ${\left[x',y'\right]}_{\gamma}$.

Let $t\in C_3\cap N_{x'},\ s\in C_3\cap N_{y'}$ and let $q$ be
a path that joins $t$ and $s$ as above. We note now that $N_{x'}$
and $N_{y'}$ are path connected, so there are paths $q_1\in N_{x'}$ and
$q_2\in N_{y'}$ joining the endpoints of $q$ with the points $x'$ and
$y'$ respectively. Clearly then the path $p=q_1\cup q\cup q_2$
joins $x'$ with $y'$ without meeting $\gamma'$.
This is however impossible, since $x'$ and $y'$
are separated by $\gamma'$. Therefore, $I$ does not cross
$\gamma$. Thus $\gamma\sm I$ is contained in a single component of
$X\sm I$, a contradiction.
\end{proof}

We return now to the proof of Lemma \ref{A.1}: 
Let $G$ be the group of homeomorphisms of $X$. 
For every $f\in G$
we have that $f(\gamma)$ separates $X$ and from the previous lemma
it follows that $f(\gamma)$ does not cross $\gamma$. Let
$S=G\cdot\gamma$. Clearly $S$ is uncountable. Let $Q$ be a
countable dense set of $X$. We define a map $R:S\to Q\times
Q\times Q$ as follows: Let $p\in S$ and $U_1,U_2,U_3$ be three
connected components of $X\sm p$. For every $U_i$ we pick an
$r_i\in Q$ and we associate $p\in S$ the triple $(r_1,r_2,r_3)$.
We remark that $R$ is $1-1$ map, which is a contradiction. 
This completes the proof of Lemma \ref{A.1}.
\end{proof}

\

\begin{figure}[h]
\begin{center}
\includegraphics[width=3in]{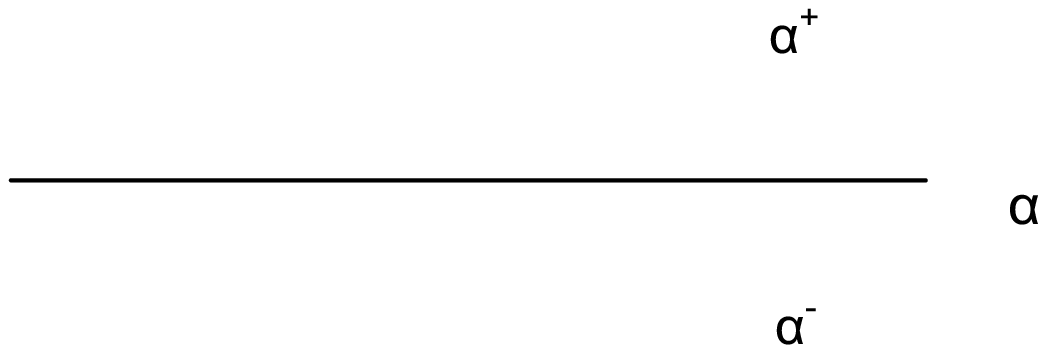}
\end{center}
\os

\label{fig:separc}
\end{figure}

Let $\gamma$ be an arc that separates $X$ in exactly two
components $C_1, C_2$ with $\partial C_1=\partial C_2=\gamma$. We
denote by $G$ the group of homeomorphisms of $X$. Let $S=G\cdot
\gamma$. It is clear that $S$ is uncountable and that every 
arc $\al\in S$ also separates $X$ in exactly two components $U_1,U_2$
such that $\partial U_1=\partial U_2=\al$. For an arc $\al\in S$ we
will denote these two components by $\al^+$ and $\al^-$ (Figure \ref{fig:separc}).

Henceforth we will consider only arcs in $S$.

\

\begin{lemma}
\label{A.2} Let $\al_1,\al_2\in S$ such that $\al_1$ crosses
$\al_2$ at $x$ (or at $[x_1,x_2]$). Then $\al_2$ crosses $\al_1$
at $x$ (or at $[x_1,x_2]$).
\end{lemma}

\begin{proof}
Suppose that there are $\al_1,\al_2\in S$ such that $\al_1$
crosses $\al_2$ at $x$ but $\al_2$ does not cross $\al_1$ at $x$.
Then there is an interval $I\subset \al_1$ containing $x$ at its
interior that lies in the closure of one of the components of
$X\smallsetminus \al_2$, say $\overline {\al_2^+}$. Clearly then
we have that $I\cap \al _2^-=\emptyset$.

Let $V\subset X$ be sufficiently small neighborhood of $x$ such
that $(\al_1\smallsetminus I)\cap V=\emptyset$. We denote by $J$
the connected component of $\al_2\cap V$ that contains the point
$x$. We pick two points $a,b\in J$ with $a<x<b$ in $\al_2$ which
are separated by $\al_1$ and let $N_a,N_b$ be connected
neighborhoods of $a$ and $b$ respectively such that
$N_a,N_b\subset V$.
As in proof of Lemma \ref{subl},
for $\ep<\min\{\diam(N_a\cap \al_2^-),\diam(N_b\cap \al_2^-)\}$,
we can find a path $p$ that
joins $a$ with $b$ without meeting $\al_1$, a contradiction. We
argue similarly if $\al_1$ crosses $\al_2$ at an interval
$[x_1,x_2]$.

\end{proof}

\

We recall now a version of Effros' Theorem (\cite {E},
\cite {H} p. 561):
\begin{theorem}
\label{Ef} For every $\ep>0$ and $x\in X$ the set $W(x,\ep)$ of
$y\in X$ such that there is a homeomorphism $h:X\to X$ with
$h(x)=y$ and $d(h(t),t)<\ep$ for all $t\in X$, is open.
\end{theorem}

\

\begin{lemma}
\label{aaaa}
There are arcs $\al=[a_1,a_2]$ and $\beta=[b_1,b_2]$ in $S$,
such that $b_1\in(a_1,a_2)_{\al}$ and if $A$ is the connected component
of $\al\cap\beta$ that contains $b_1$, then $a_1,a_2\not\in A$.
\end{lemma}

\begin{proof}
We will need the following:
\begin{sub}
\label{cross}
Let $\al\in S$. Then there is an arc $\beta\in S$ that crosses $\al$.
\end{sub}

\begin{proof}

Let  $\al, \gamma \in S,\ c\in\partial\gamma,\
a\in\al\sm\partial\al$ and $g\in G$ such that $gc=a$. By the
definition of $S$ it is not possible that $g\gamma \subset \al
$, since $\al$ separates $X$ in exactly two connected components.
Assume now that $\alpha$ does not cross $g\gamma$.
We denote by $A$ the connected component of $\alpha\cap g\gamma$
that contains $a$ and let $\partial \alpha=\{a_1,a_2\}$.

We distinguish two cases:
Suppose that $a_1,a_2\not\in A$. Let $z\in g\gamma$ such that
$(z,gc)_{g\gamma}$ lies in the closure of one of the components of
$X\sm\al$, say $\al^+$. Let $z'\in(z,gc)_{g\gamma}\sm\al$ and
$\ep>0$ such that $B_{\ep}(z')\subset \al^+$.
By Theorem \ref{Ef} there is a $\delta>0$ such that
$B_{\delta}(a)\subset W(a,\ep)$. Let $y\in B_{\delta}(a)\cap
\al^-$ (Figure \ref{befef1}).

\begin{figure}[h]
\begin{center}
\includegraphics[width=3.5in]{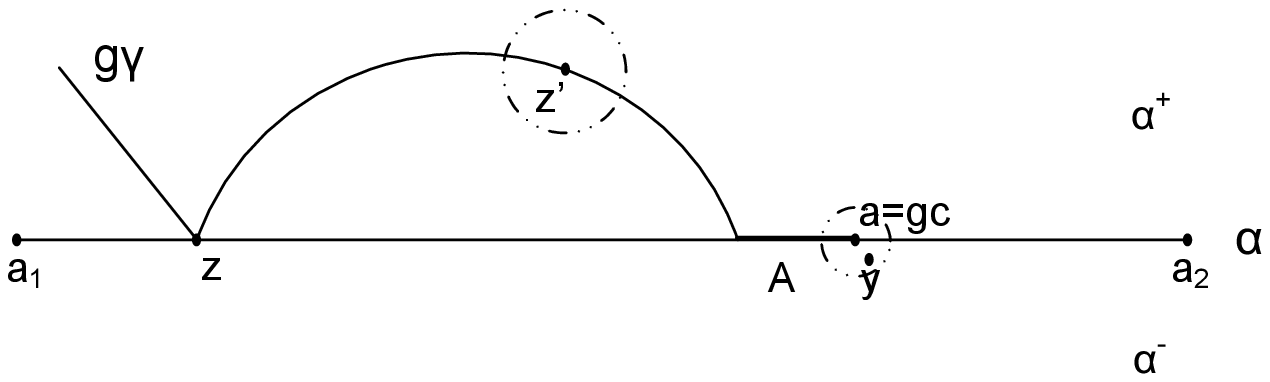}
\end{center}
\os
\label{befef1}
\end{figure}

Then there is a homeomorphism, $h\in G$, with $h(a)=y$ such that
$d(t,h(t))<\ep$ for every $t\in X$. We consider the arc
$\beta=h(g\gamma)$. Then clearly $\beta$ crosses $\al$, since
$h(z')\in \al^+$ and $h(a)\in \al^-$.

Suppose now that $a_2\in A$. We consider the homeomorphism $h\in G$ of
the previous case. If $a_2\not\in h(A)$, then clearly  Lemma \ref{cross} is proved. So
let $a_2\in h(A)$ and $\ep'<\min\{\ep, \frac{1}{2}d(\al,h(gc))\}$.
 As before, by Theorem \ref{Ef}, there is a $\delta'>0$ such that
$B_{\delta'}(a_2)\subset W(a_2,\ep')$.
Let $y'\in B_{\delta'}(a_2)\cap\al^+$ (Figure \ref{befef2}).
Then there is an $h'\in G$ with $h'(a_2)=y'$ and
$d(t,h'(t))<\ep'$ for every $t\in X$.
It is obvious now that $\al$ crosses $h'(\beta)$.

\begin{figure}[h]
\begin{center}
\includegraphics[width=3in]{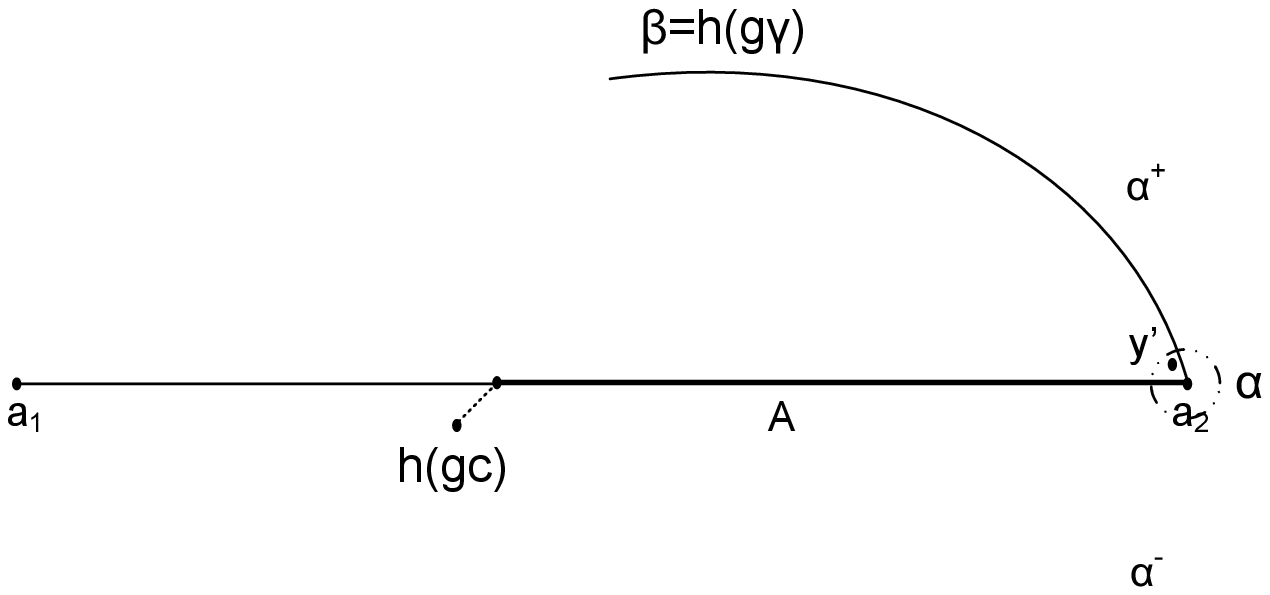}
\end{center}
\os
\label{befef2}
\end{figure}
\end{proof}

Let $\al=[a_1,a_2]\in S$.
By Lemma \ref{cross} there is an arc $\beta=[b_1,b_2]$
that crosses $\al$ at $x\in\al\cap\beta$.
Without loss of generality, suppose that $x$
is the endpoint of a cross interval $I$ of $\al,\beta$.
Let $\gamma\in S,\ c\in\partial \gamma$ and $g\in G$ such that $gc=x$.
We denote by $A$ the connected component of $g\gamma\cap\al$
that contains $c$ and similarly by $B$ the component of
$g\gamma\cap\beta$ that contains $c$.
Clearly if $a_1,a_2\notin A$, then Lemma \ref{aaaa} is proved.
Otherwise, we note that if $A$ contains one of the endpoints of $\al$,
then $b_1,b_2\notin B$, since $I$ is a cross interval of $\al,\beta$.
So in this case, the required arcs are $g\gamma$ and $\beta$.
\end{proof}

\

\noindent We return to the proof of Theorem \ref{continua}.

\

Let $\al=[a_1,a_2],\beta=[b_1,b_2]\in S$ be paths as in Lemma \ref{aaaa},
that is $b_1\in(a_1,a_2)_{\al}$ and if $A$ is the connected component
of $\al\cap\beta$ that contains $b_1$, then $a_1,a_2\notin A$.
Let $t_1,t_2\in \al\sm \beta$ such that $b_1\in {(t_1,t_2)}_{\al}$ and let
$p_1,p_2$ be paths joining $t_1,t_2$ in $\overline{ \al^+}$ and
$\overline{ \al^-}$ respectively (the points $t_i$ exist since $a_1,a_2\notin A$).
We pick $p_i$ such that $p_i\cap
\al$ has exactly two connected components neither of which
intersects $\beta$ (this can be achieved using Lemma \ref{A.0'}
for $\ep<\frac{1}{2}\min\{d(t_1,\partial A_1),d(t_2,\partial
A_2)\}$, where $A_i$ is the connected component of $\al\sm \beta$
that contains $t_i$ and $\partial A_i$ is its boundary in $\al$).

Let $\ep>0$ with $\ep<\frac{1}{2}d(A,p_1\cup p_2)$. As in proof of
Lemma \ref{cross}, using Theorem \ref{Ef}, we can find a
homeomorphism $h\in G$ such that $h(\beta)$ crosses $\al$ at
$x\in\al\cap h(\beta)$, with $d(A,x)<\ep$. Then we remark that
$x\in {(t_1,t_2)}_{\alpha }$ and that the subarc of $h\beta$ with
endpoints $x$ and $hb_1$ does not intersect $p_1\cup p_2$.

\begin{figure}[h]
\begin{center}
\includegraphics[width=3.5in]{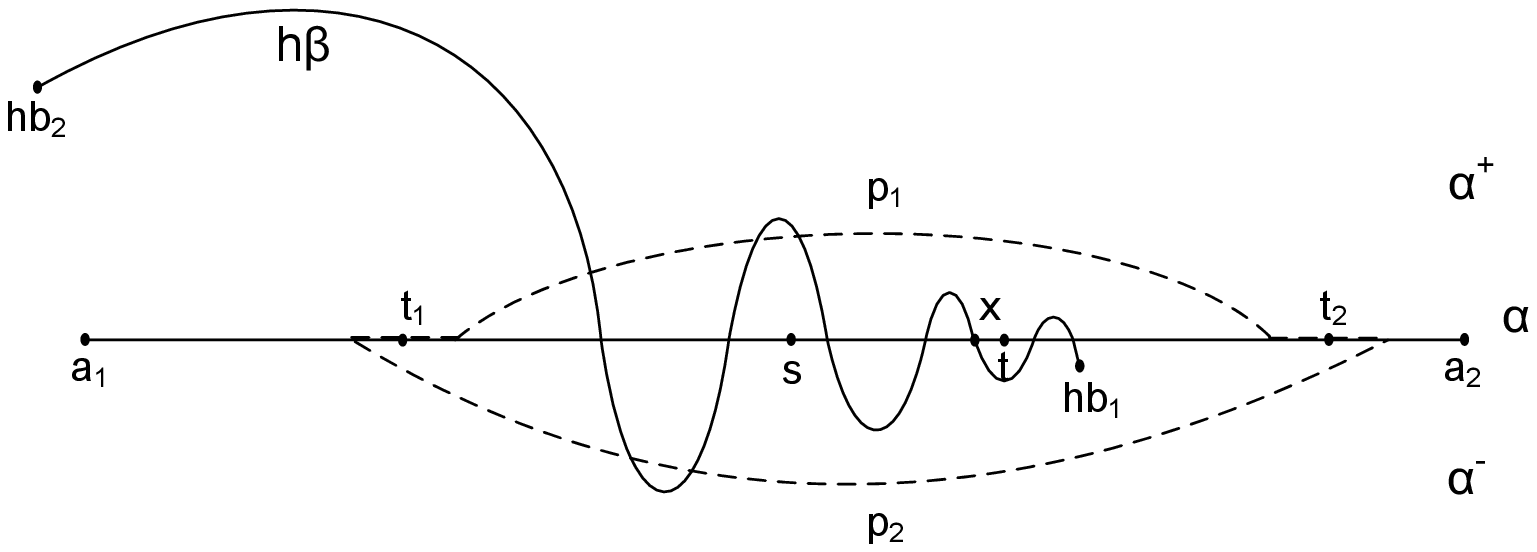}
\end{center}
\os
\label{aftef}
\end{figure}

We pick now points $s\in {(t_1,x)}_{\al}\sm (p_1\cup p_2)$ and
$t\in {(x,t_2)}_{\al}\sm (p_1\cup p_2)$ which are
separated by $h\beta$ so that they satisfy the following:
If $y$ is a cross point of $\al$ and $h\beta$, lying in ${[s,t]}_{\al}$,
then the subarc ${\left[y,hb_1\right]}_{h\beta}$ does not intersect the
paths $p_1,p_2$. Such points exist by definition of $x$ (Figure \ref{aftef}).

We consider now the closed paths $p_1\cup [t_1,t_2]_{\alpha }$ and
$p_2\cup [t_1,t_2]_{\alpha }$. Let $D_1,D_2$ be discs and let $f_1:D_1\to\overline{\al ^+}$,
$f_2:D_2\to \overline{\al^-}$ be maps so that $f_1(\partial D_1)=p_1\cup [t_1,t_2]_{\alpha }$,
$f_2(\partial D_2)=p_2\cup [t_1,t_2]_{\alpha }$ (such maps exist,
since $\overline{\al^+}$ and $\overline{\al^-}$ are simply connected by Lemma \ref{simplycon}).

We `glue' $D_1,D_2$ along $[t_1,t_2]_{\alpha }$ and we obtain a
disc $D$ and a map $f:D\to X$ with $f(\partial D)=p_1\cup p_2$.
More precisely, we consider the disc $D=D_1\sqcup D_2/\thicksim$,
where $\thicksim$ is defined as follows:
$x_1\thicksim x_2$ if and only if $x_1\in\partial D_1,\ x_2\in\partial D_2$
and $f_1(x_1)=f_2(x_2)$. Finally, we define $f:D\to X$ as:

\[f(t)=\left\{\begin{array}{ll}
f_1(t), & \mbox{if $t \in D_1$}, \\
f_2(t), & \mbox{if $t \in D_2$}.
\end{array}\right.\]

\begin{figure}[h]
\begin{center}
\includegraphics[width=2.7in]{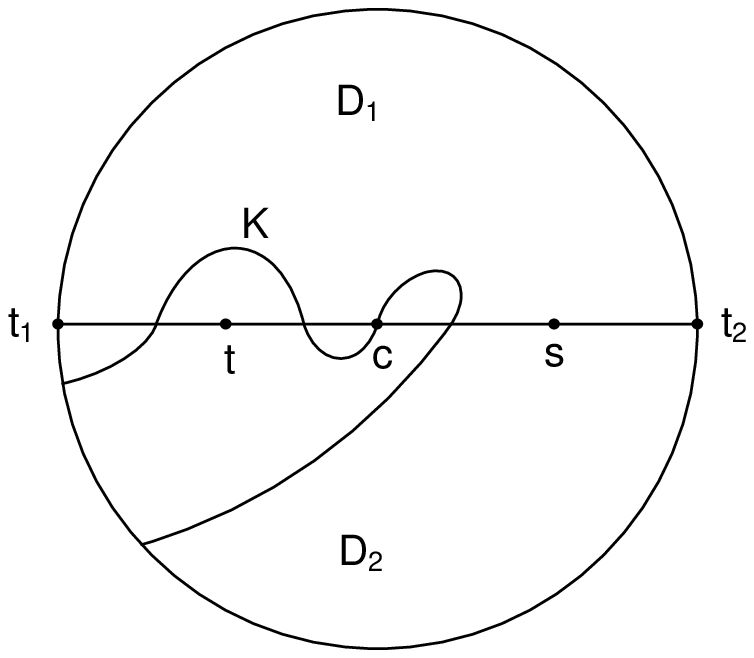}
\end{center}
\os
\label{diskos}
\end{figure}

By abuse of notation we identify points of $[t_1,t_2]_{\alpha }$
in $D$ with their image under $f$.
We note that the interior, say $U$, of $D$ is homeomorphic to
$\mathbb{R}^2$ and since $t,s$ are separated by
$h\beta$ in $X$, it follows by Alexander's lemma that $t,s$ are
separated in $U$ by a connected component of $f^{-1}(h\beta)\cap
U$. We call this component $K$ (Figure \ref{diskos}).

Clearly $f(K)$ is a subarc of $h\beta$ that contains cross points
or cross intervals of $h\beta$ with $[s,t]_{\alpha }$. Let $c$ be
such a cross point. Then we can write $f(K)$ as $f(K)=I_1\cup
I_2$, where $I_i,\, i=1,2$, are (connected) subarcs of $h\beta$,
such that $I_1\cap I_2=c$. Furthermore, at least one of $I_1,I_2$
does not intersect $p_1\cup p_2$ (this is by our choice of $h$ and
$c$). It follows that at least one of $f^{-1}(I_1)\cap U,
f^{-1}(I_2)\cap U$ is compact.

We set $I_1'=I_1\sm c,\, I_2'=I_2\sm c$.
We will define two sets $K_1,K_2$ such that the following are
satisfied: $K_1,K_2$ are closed subsets of $U$ that contain $f^{-1}(I_1')$
and $f^{-1}(I_2')$ respectively, $K_1\cap K_2$ is connected contained in  $f^{-1}(c)$ and $K_1\cup K_2=K$.

\

We consider the connected components of $f^{-1}(c)\cap K$. We
remark that there is exactly one component of $f^{-1}(c)\cap K$,
say $C$, that intersects both $D_1$ and $D_2$.

Let now $C_1$ be a connected component of $f^{-1}(c)\cap K$ different
from $C$ and suppose that $C_1\subset D_1$. We consider the closure
of $K$, $\overline K$, in $D_1\cup D_2$. $\overline K$ is connected
 thus the closure of the component of $\overline K\cap D_1$ containing
$C_1$ intersects $\partial D_1$. Indeed, we consider the set
$\overline{K}\cap  (D_1-[t_1,t_2]_{\alpha})$
 as an open subset
of the continuum $\overline{K}$. Let $K'$ be the component of
$\overline {K}\cap (D_1-[t_1,t_2]_{\alpha})$ that contains $C_1$.
We recall that if $U$ is an open subset of a continuum and $C$ is
a component of $U$ then the frontier of $U$ contains a limit point
of $C$ (Theorem 2.16, p.47 of \cite{H-Y}). It follows that the
closure of $K'$ intersects $[t_1,t_2]_{\alpha}$.

Therefore, we have that
$f(K')\subset \overline{\al^+}$ so $f(K')\subset I_1$ or
$f(K')\subset I_2$. We remark that if $f(K')\subset I_1$ then a
non trivial interval of $I_1$ containing $c$ lies in
$\overline{\al^+}$.

We have a similar conclusion if $f(K')\subset I_2$. Therefore if a
connected component of $f^{-1}(c)\cap K$ different from $C$ lying
in $D_1$ intersects the closure of both
$f^{-1}(I_1'),f^{-1}(I_2')$ we have that an open interval of $I_1$
around $c$ lies in $\overline{\al^+}$. This is impossible since
$c$ is a cross point.
We argue similarly for connected components of $f^{-1}(c)\cap K$
contained in $D_2$.

We conclude that the union of the components of
$f^{-1}(c)\cap K$ which lie in $D_1$, intersect exactly one of
$f^{-1}(I_1'),f^{-1}(I_2')$.
Clearly the same is true for the union of the components of $f^{-1}(c)\cap K$
contained in $D_2$.
In particular exactly one of the following two holds:

\

1. If $C_1$ is a connected component of $f^{-1}(c)\cap K$
different from $C$ lying in $D_1$ then the component of $D_1\cap
K$ containing $C_1$ intersects $f^{-1}(I_1')$, while if $C_1$ lies
in $D_2$  the component of $D_2\cap K$ containing $C_1$ intersects
$f^{-1}(I_2')$.

2. If $C_1$ is a connected component of $f^{-1}(c)\cap K$
different from $C$ lying in $D_1$ then the component of $D_1\cap
K$ containing $C_1$ intersects $f^{-1}(I_2')$, while if $C_1$ lies
in $D_2$  the component of $D_2\cap K$ containing $C_1$ intersects
$f^{-1}(I_1')$.

\

Assume that we are in the first case. Then we define $K_1$ to be
the union of the components of $f^{-1}(c)\cap K$ intersecting
$D_1$ together with $f^{-1}(I_1')$. We define $K_2$ to be the
union of the components of $f^{-1}(c)\cap K$ intersecting $D_2$
together with $f^{-1}(I_2')$. It is clear that $K_1,K_2$ are
closed and that $K_1\cap K_2=C$, $K_1\cup K_2=K$. Since $K$ is
connected, $K_1,K_2$ are connected too.
We define $K_1,K_2$ similarly in the second case.

We note now that at least one of $K_1,K_2$ is compact subset of
$U$, thus bounded in $U$. Since $K$ separates $s,t$ and $K_1,K_2$
are closed subsets of $D$, applying Alexander's lemma for the
plane we have that at least one of $K_1,K_2$ separates $s,t$ in
$U$.

It follows that either $f^{-1}(I_1)$ or $f^{-1}(I_2)$ separates
$s,t$. We remark that the same argument holds in the case $c$ is
replaced by a cross interval $J$: We have then that $I=I_1\cup
I_2$ with $I_1\cap I_2=J$ and as before either $f^{-1}(I_1)$ or
$f^{-1}(I_2)$ separate $s,t$ in $U$. Now we can continue
subdividing intervals along cross points (cross intervals) that
lie in $[s,t]_{\alpha }$ as follows: Let's say that $f^{-1}(I_1)$
separates $s,t$. We have that there is a connected component of
$f^{-1}(I_1)$, say $M$, that separates them. We note that $f(M)$
is a subinterval of $I_1$ and if there is a cross point (or cross
interval) of $[t,s]_{\al },h\beta$ contained in its interior, we
repeat the previous procedure replacing $K$ by $M$. If not we have
a contradiction. Therefore, either $s,t$ are separated in $U$ by
the inverse image of an interval $f(K)$ of $h\beta$ which does not
contain in its interior any cross point of $h\beta, \alpha $ lying
in $[s,t]_{\alpha }$, or by iterating this procedure we conclude
that the inverse images under $f$ of intervals of $h\beta $ of
arbitrarily small diameter separate $s$ from $t$ in $U$. It is
clear that both are impossible, so the theorem is proven.

\end{document}